\documentclass[12pt,reqno]{amsart}
\usepackage{amssymb,amsmath}
\def\la{\lambda}

\def\la{\lambda}

\def\g{\gamma}
\def\d{\delta}

\def\({\left(}

\def\){\right)}

\def\eb{\varepsilon}

\def\eb{\varepsilon}

\def\la{\lambda}

\def\R {\mathbb{R}}

\def\B {{\mathcal B}}

\def\E {{\mathcal E}}

\def\E {{\mathcal E}}

\def\W {{\mathcal W}}

\def\de{\delta}

\def \ga{\gamma}
\def \and{\qquad\text{and}\qquad}

\newtheorem{proposition}{Proposition}[section]
\newtheorem{theorem}[proposition]{Theorem}

\newtheorem{lemma}[proposition]{Lemma}
\theoremstyle{definition}

\newtheorem{remark}[proposition]{Remark}

\numberwithin{equation}{section}

\def \no#1#2#3 {{\bf #1} (#3), #2.}
\def \eds#1#2#3 {#1, #2, #3.}
\def\be{\begin{equation}}
\def\ee{\end{equation}}

\title[Stabilization of Marine Riser Equation] {Stabilization of  
  the Marine Riser model by  controllers depending on finitely many parameters}
\author[] {V.K. Kalantarov,   A.A. Namazov and E. S. Titi}

\address{(V.K.Kalantarov) Department of Mathematics,
\newline\indent Ko{\c c} University,  Istanbul, Türkiye;
\newline\indent  Department of Engineering Mathematics and Artificial Intelligence,  
\newline \indent Azerbaijan Technical University, Baku, Azerbaijan. \ e-mail: vkalantarov@ku.edu.tr}
\address{(A.A. Namazov)  Department of Engineering Mathematics and Artificial Intelligence, 
\newline\indent Azerbaijan Technical University, Baku, Azerbaijan.  \ e-mail: atif.namazov@aztu.edu.az}

\address{{E.S. Titi} Department of Mathematics, Texas A M University, College Station,\newline\indent  TX 77843, USA;
\newline\indent  Department of Applied Mathematics and Theoretical Physics, University of
\newline\indent  Cambridge, Cambridge CB3 0WA, UK;
\newline\indent  Department of Computer Science and Applied Mathematics, Weizmann Institute of\newline\indent  Science, Rehovot 7610001, Israel.  \ e-mail: edriss.titi@maths.cam.ac.uk}

\date{ 27 February 2026}
 \keywords{}
\usepackage[dvips]{color}
\begin{document}

\begin{abstract} 
We prove global stabilization of the marine riser models using a feedback controller that depend on finitely many finite-volume elements 
and finitely many nodal observables. Our approach is based on a  feedback control design for dissipative nonlinear
 partial differential equations, inspired by the methodology introduced in [Evol. Equ. Control Theory, {\bf 3} (2014), 579--594]. 
The proposed control strategy ensures asymptotic stabilization while maintaining computational feasibility, making it suitable for practical applications.

\end{abstract}

 \maketitle

\begin{center}

\textit{This paper is dedicated to Professor Roger Temam, on the occasion of his 85th birthday.}

\end{center}

\section{Introduction}

We study the problem of stabilization to the zero stationary state of the marine riser equation
\begin{equation}\label{C}
mu_{tt} + ku_{xxxx} - \left[ {a(x)u_x } \right]_x + \gamma u_{tx}
+bu_t \left| {u_t } \right|^p = 0,\ x\in (0,L), t>0,
\end{equation}
subject to the homogeneous joint Dirichlet and Neumann  boundary conditions
\begin{equation}\label{1.5}
u(0,t)=u_{x}(0,t)=u(L,t)=u_{x}(L,t)=0, \quad t>0,
\end{equation}
where $k>0$ is the  flexural   rigidity of the riser,  $b>0$ is the nonlinear drag force and $m>0$ represents the mass
line density, $\gamma $ is the Coriolis force parameter,  $a(\cdot)\in  C^1[0, L]$
 describes the effective tension coefficient and $u| {u_t } |^p, \ p\in [1,\infty)$ is the nonlinear damping term.

The nonlinear wave equation  \eqref{C} wıth $p=1$ arises in the mathematical model of marine risers - long, slender pipes used in offshore drilling to transport oil or gas from the seabed to surface platforms. Such structures are subject to complex interactions with surrounding fluids and external forces, including gravity, ocean currents, and drilling-induced vibrations. Their safe operation critically depends on stability: excessive vibrations or negative effective tension may lead to buckling, compromising the integrity of the riser. The study of the dynamical behavior of equation \eqref{C} is therefore not only of mathematical interest, but also of direct importance in offshore engineering. At the same time, it fits naturally into the broader mathematical framework of feedback stabilization of nonlinear dissipative partial differential equations, where one seeks to design effective control mechanisms to suppress instabilities and ensure long-time stability.\\

Earlier studies of the marine riser model (see e.g. \cite{Ko,KaKu}) established stability of the zero solution by Lyapunov methods when the effective tension, $a(x)$, is positive, and provided decay in time  estimates (polynomial or exponential, depending on the damping).  Specifically, in \cite{Ko}   the Lyapunov function technique  is used to show that
 the zero solution of  equation \eqref{C}  under the homogeneous Dirichlet's boundary conditions
\be\label{Dir}
u(0,t)=u_{xx}(0,t)=u(L,t)=u_{xx}(L,t)=0, t>0,
\ee
 is stable when the coefficient of effective tension is positive.
The polynomial decay in time estimate for solutions of this problem is
established \cite{KaKu} when the effective tension
 $a(x)$ is
a $C^1 \left[ {0,L} \right]$ function that satisfies the conditions  
\be\label{ax} -a_0\le  a(x) \le a_1,   \
\ \forall x \in[0,L], \quad \hbox{with} \quad  a_0 \ge0, \  a_1>0 ,
\ee 
and
\be\label{d0}
 d_0:=k-a_0\frac{L^2}{\pi^2}>0.
\ee
Similar  results for the multidimensional marine riser equation are obtained in \cite{KaKu} and are adapted in  \cite{Aa} and \cite{Gur}. Furthermore, the problem of structural stability of marine riser equation is discussed in  \cite{CKU1}.\\

 However, it is well known that if the effective tension becomes negative, the riser may buckle and stability is lost.
Thus a natural question arises: is it possible to stabilize the system employing some feedback controllers. The problem of   feedback stabilization of nonlinear evolutionary PDEs has been extensively investigated in the literature  (see, e.g. \cite{Coron, Fur, Rus, Trig, Zua} and references therein.) In particular, the development of finite-dimensional feedback controllers (see e.g. \cite{AzTi,KaTi1,KaTi,LuTi}) has attracted significant attention in recent years, as such controllers are both mathematically efficient and practically implementable in engineering applications.\\

The design of finite-dimensional controllers for PDEs has a long history, going back to works such as \cite{Ba,Trig}, and has been extensively developed in the last decades for fluid flows, wave equations, and reaction–diffusion systems (see, e.g., \cite{BaTa, Barbu, Cheb,Munt, Sel, Yan} and references therein). Of particular importance for our approach is the feedback control framework introduced in \cite{AzTi}, which proposed an efficient finite-dimensional algorithm for globally stabilizing dissipative PDEs. This method has since been applied to a wide class of systems, including nonlinear wave equations, the Boussinesq, the Kuramoto–Sivashinsky equations, the complex Ginzburg–Landau equation, Navier–Stokes–Voigt models, and chevron pattern equations (see \cite{AzTi,KaTi,KaTi1,KaOz,KKV1,LuTi} and references therein).\\

The Azouani-Titi \cite{AzTi} approach is closely related to the continuous data assimilation algorithm of Azouani, Olson, and Titi \cite{AOT}, in which feedback terms are introduced directly at the PDE level to incorporate partial observational data. This perspective connects control theory with state estimation and model synchronization, and has inspired numerous studies in both theoretical and applied contexts (\cite{Far,FJT,Mar,OlTi,Lar1,LhPr}). 

Recently, similar ideas have been applied to the stabilization of the marine riser model with {negative effective tension}. In particular, \cite{AKN} proved the global asymptotic stabilization of the zero stationary state using feedback controllers depending on finitely many Fourier modes. This result demonstrates that even when the uncontrolled system is unstable, it can be stabilized through appropriately designed low-dimensional control laws.

In this paper, we extend the results of \cite{AKN} and develop new feedback control strategies for the marine riser equation \eqref{C}–\eqref{1.5}. Our main contributions can be summarized as follows:

\begin{itemize}
  \item  We establish {global asymptotic stabilization} of the marine riser equation using feedback controllers depending on finitely many {finite-volume elements} rather than Fourier modes. These feedbacks are more natural from a practical experimental and computational viewpoint. {From an applied perspective}, the proposed feedback design is computationally efficient and suitable for real-time implementation, as it requires only coarse measurements of the system state. It, therefore, offers a promising approach to vibration control in flexible offshore structures.
  \item  We provide explicit conditions on the controller parameters ensuring stabilization and derive uniform energy estimates for the closed-loop system.
  \item  We analyze the influence of both linear and nonlinear damping mechanisms on the decay rates of the energy, identifying parameter regimes leading to polynomial or exponential in time stabilization rate.
  \item  Our analysis combines Lyapunov functional techniques, compactness arguments, and the finite-dimensional feedback framework introduced in \cite{AzTi}, adapted here to a fourth-order wave-type equation with variable coefficients. The results obtained contribute to the ongoing development of control and stabilization theory for nonlinear dissipative nonlinear wave equations, particularly for models involving both tension and bending effects.
\end{itemize}

The remainder of the paper is organized as follows. In Section~2 we present the functional setting, main assumptions, and preliminary lemmas. Section~3 is devoted to the design of the feedback controllers and the proof of global stabilization results. Here  we study the influence of nonlinear damping and provide explicit decay rate estimates. Finally, in Section~4   feedback controllers  based on finite-volume elements are utilized to stabilize the zero stationary state of the linear equation 
$$
mu_{tt} + ku_{xxxx} - \left[ {a(x)u_x } \right]_x + \gamma u_{tx}
+bu_t   = 0, \ \ x\in (0,L), t>0,
$$
modeling  dynamics of the  marine riser conveying
fluid under the homogeneous Dirichlet's  boundary conditions \eqref{1.5}.

\section{Preliminaries}
\indent In what follows we use  the Young's inequality with $\eb>0$
\be\label{Cauchy}
| XY| \le \frac{\eb}2X^2+\frac1{2\eb}Y^2, \ \ \forall X,Y \in \R,
\ee
the  Poincar{\'e} inequality
\be\label{Po}
\|v\|^2\leq \lambda_1^{-1}\| v_x\|^2, \ \ \forall v\in H_0^1(0,L),
\ee
where $\lambda_1 >0$ is the first eigenvalue of the operator $-\partial_{xx}$ on $(0,L)$ subject to Dirichlet homogenous boundary condition, and the Sobolev inequality
\be\label{Sob}
\|v\|_{L^q(0,L)}\le C_L\|v_x\|, \ q\in[1,\infty], \ \forall v\in H_0^1(0,L).
\ee
Here and in the sequal  $(\cdot, \cdot)$ and $\|\cdot\|$  denote the  inner product and the norm of $L^2(0,L)$, respectively. 

Moreover, we will be using  following lemma in the proof of our main results
\begin{lemma}\label{pr}(see \cite{AzTi})  Let $\phi\in H^1(0,L)$.  Then, for any positive integer $N$, with $h:=\frac LN$, one has
\be\label{p1}
\|\phi-\sum\limits_{k=1}^N\overline{\phi}_k\chi_{J_k}(\cdot)\|\leq h
\|\phi_x\|, \ee and 
\be\label{p2} \|\phi\|^2\leq
h\sum\limits_{k=1}^N\overline{\phi}_k^2+h^2\|\phi_x\|^2,
\ee 
where $J_k:=\left[(k-1)\frac LN, k\frac LN\right),$ for $k=1,2,\cdots N-1$ and $J_N=[\frac{N-1}{N}L, L]$,

$$\overline{\phi}_k:=\frac1{|J_k|}\int\limits_{J_k}\phi(x)dx,
$$ 
and
$\chi_{J_k}(x)$ is the characteristic function of the interval
$J_k$. \\
\end{lemma}

\section {Stabilization of the equation with nonlinear damping term} In this section we employ  feedback controllers  based on finite-volume elements to stabilize the marine riser equation, \eqref{C}–\eqref{1.5}, with nonlinear damping term. That is, 
 we consider the following closed-loop feedback control problem
\begin{equation}\label{kg1}
mu_{tt}  + ku_{xxxx}  - \left[ {a(x)u_x } \right]_x  + \g u_{tx}
+ bu_t \left| {u_t } \right|^p=
-\mu\sum\limits_{k=1}^N\overline{u}_k(t)\chi_{J_k}(x),   \ x \in (0,L), \ t>0, 
\end{equation}
\begin{equation}\label{kg2}
{u}(0,t)={u}(L,t)=u_{x}(0,t)=u_{x}(L,t)=0, \ t>0,
\end{equation}
\begin{equation}\label{kg3}
u(x,0)=u_0(x), \ \ u_t (x,0)=u_1(x), \ \ x\in (0,L),
\end{equation}
where $\mu>0$, is the feedback nudging parameter, and $a(\cdot)\in C^1[0,L]$ satisfies only  condition \eqref{ax}.

First let us note that by using the standard Faedo-Galerkin method one can prove the following theorem about global unique solvability of problem \eqref{kg1}-\eqref{kg3}  (see, e.g,   \cite{ Lions}, \cite{LiSt},\cite{Mar})
\begin{theorem} \label{exun} Let $u_0\in H_0^2(0,L), \ \ u_1\in  H_0^2(0;L)$ and $T>0$. Then  problem \eqref{kg1}-\eqref{kg3} has a unique  solution 
$$
u\in L^\infty(0,T;H_0^2(0;L)), \ u_t\in L^\infty(0,T;  H^2_0(0,L)) ,
$$
that satisfies equation \eqref{kg1}  in the sense of distributions and attains the initial conditions in $L^2$ sense. Moreover, the following estimates are satisfied
\be\label{Esta0}
\|u_t(t)\|^2+\|u_{xx}(t)\|^2+\|u(t)\|_{L^{p+2}}^{p+2}(Q_T)\le C_0,\
\ee
\be\label{Esta1}
\|u_{tx}(t)\|^2+\|u_{xxx}(t)\|^2\le C_1,\
\ee

where $Q_T:=(0,L)\times (0,T)$  and  $C_0, C_1$ are  positive constants independent of $T.$
\end{theorem}

More details about existence, uniqueness and regularity  of solutions to initial boundary value problems for  wave type equations, with linear and nonlinear damping,  can be found, e.g., in  \cite{ChLa}, \cite{Tem} and in  references therein. \\
To prove  estimates \eqref{Esta0} and \eqref{Esta1},  as well as other estimates of solutions to problem \eqref{kg1}-\eqref{kg3} we  restrict ourselves to formal derivation of the estimates. However, these estimates can be rigorously justified by establishing them first to, e.g.,  the Faedo-Galerkin  approximation system and then passing to the limit after using the appropriate compactness theorems.\\

First, we  derive the following energy equality by taking the inner product of 
 the equation  \eqref{kg1} with $u_t$ :
\begin{multline}\label{E1}
\frac d{dt}\left[\frac m2 \|u_t(t)\|^2+\frac k2\|u_{xx}(t)\|^2
+\frac12\int_0^La(x)u^2_x(x,t)dx+\frac{\mu h}2\B_N(u(t))\right]\\
=-b\int_0^L|u_t(x,t)|^{p+2}dx.
\end{multline}
Here, and in what follows we denote 
\be\label{BNu}\B_N(u(t)):=\sum\limits_{k=1}^N\overline{u}_k^2(t), \ \ \overline{u}_k(t)= \frac1{h}\int\limits_{J_k}u(x,t)dx. \ee
Next, we take the inner product of \eqref{kg1} with $u$ to obtain

\begin{multline}\label{E2}
m\frac d{dt}(u(t),u_t(t))-m\|u_t(t)\|^2+k\|u_{xx}(t)\|^2+\int_0^La(x)u^2_x(x,t)dx\\
-\gamma(u_x(t),u_t(t))+b\int_0^Lu(x,t)u_t(x,t)|u_t(x,t)|^pdx+\mu h\B_N(u(t))=0.
\end{multline}

Then we multiply \eqref{E2} by some $\de>0$  (to be choosen below) and add the resultant to \eqref{E1} to get
\begin{multline}\label{dtE}
\frac d{dt} E(u(t))-\de m\|u_t(t)\|^2+\de k\|u_{xx}(t)\|^2+\de\int_0^La(x)u_x^2(x,t)dx-\de\ga(u_x(t),u_t(t))\\+b\de\int_0^Lu(x,t)u_t(x,t)|u_t(x,t)|^pdx
-b\int_0^L|u_t(x,t)|^{p+2}dx+\de\mu h\B_N(u(t))=0,
\end{multline}
where
$$
 E(u(t)):=\frac m2\|u_t(t)\|^2+\frac k2\|u_{xx}(t)\|^2+\frac12\int_0^La(x)u^2_x(x,t)dx+\frac{\mu h}2\B_N(u(t))+\de m(u(t),u_t(t)).
$$
By virtue of   the following interpolation inequality
\be\label{int1}
\|u_x\|^2\le \|u\|\|u_{xx}\|,
\ee 
which is valid for each function $u\in H^2(0,L)\cap H_0^1(0,L)$, and Poincar{\'e} inequality, \eqref{Po}, one has 
\be\label{Po1x}
\|u\|^2\le \la^{-2}_1\|u_{xx}\|^2.
\ee
Using \eqref{Po1x} together with  Cauchy-Schwarz and Young's inequalities imply
\be\label{uut}
\d m|(u(t), u_t(t))|\le \d m\|u_t(t)\| \|u(t)\| \le \frac m4\|u_t(t)\|^2 + \frac{m\d^2}{\la^2_1}\|u_{xx}(t)\|^2.
\ee
Moreover, from condition \eqref{ax}, together with \eqref{int1}, \eqref{p2} and \eqref{Po1x} we also  have
\begin{multline}\label{uxao}
-\frac12\int_0^La(x)u_x^2(x,t)dx\le\frac{a_0}2\|u_x(t)\|^2\le \frac{a_0}2 \|u(t)\|\|u_{xx}(t)\|\le \frac k4\|u_{xx}(t)\|^2+\frac{a^2_0}{4k}\|u(t)\|^2\\
\le
\frac k4\|u_{xx}(t)\|^2+\frac{a^2_0}{4k}h\B_N(u(t))+
\frac{h^2a_0^2}{4k\la_1}\|u_{xx}(t)\|^2.
\end{multline}
Putting things together imply the following lower estimate  for $E(t)$ 

\be\label{En1}
E(u(t))\ge \frac m4\|u_t(t)\|^2+\left[\frac k4- \frac{h^2a^2_0}{4k\la_1} -
\frac{m\d^2}{\la^2_1} \right]\|u_{xx}(t)\|^2+\left[\frac{\mu h}2-\frac{a^2_0}{4k}h\right]\B_N(u(t)).
\ee
Here we choose 
\be\label{del1}\d \in \left(0, \frac{\la_1\sqrt{k}}{4\sqrt{m}}\right).\ee 
Consequently,  if we assume that $h$ is small enough and $\mu$ is large enough such that
\be\label{hmu1}
h\le \frac{k\sqrt{\la_1}}{2a_0}, \ \ \mbox{and} \ \ \mu\ge\frac{a_0^2}{2k},
\ee
then one has
\be\label{Etge}
E(u(t))\ge \frac m4\|u_t(t)\|^2+\frac k8\|u_{xx}(t)\|^2.
\ee
Using similar arguments one can also establish the following  inequality
$$
\de\ga|(u_x(t),u_t(t))|\le \frac12\|u_t(t)\|^2+\frac12\de^2\ga^2\la_1^{-1}\|u_{xx}(t)\|^2.
$$
Using the above we deduce from  \eqref{dtE}   the inequality
 \begin{multline}\label{dtE9}
\frac d{dt} E(u(t))-(\de m+\frac12)\|u_t(t)\|^2+\left(\de k-\frac12\de^2\ga^2\la_1^{-1}\right)
\|u_{xx}(t)\|^2+\de\int_0^La(x)u_x^2(x,t)dx\\+b\de\int_0^Lu(x
,t)u_t(x,t)|u_t(x,t)|^pdx
+b\int_0^L|u_t(x,t)|^{p+2}dx+\de\mu h\B_N(u(t))\le0.
\end{multline}
Suppose in addition to \eqref{del1} that 
\be\label{del}
0<\de<\frac{2k\la_1}{\ga^2}.
\ee
Then we obtain from  \eqref{dtE9} the following inequality
\begin{multline}\label{dtE9a}
\frac d{dt} E(u(t))\le (\de m+\frac12)\|u_t(t)\|^2+m\|u_t(t)\|^2\\
-m\|u_t(t)\|^2-\de D_0\|u_{xx}(t)\|^2-\de \int_0^La(x)u^2_x(x,t)dx
-\de\mu h\B_N(u(t))\\-b\de\int_0^Lu(x,t)u_t(x,t)|u_t(x,t)|^pdx
-b\int_0^L|u_t(x,t)|^{p+2}dx,
\end{multline}
where 
  $$D_0:= k-\frac12\de\ga^2\la_1^{-1}.$$

Consider the function
$$
\E(u(t)):=\frac m2\|u_t(t)\|^2+\frac k2\|u_{xx}(t)\|^2+\frac12\int_0^La(x)u^2_x(x,t)dx
+\frac {\mu h}2\B_N(u(t))).
$$
Since $a(x)\ge-a_0$, thanks to \eqref{ax}, then by using \eqref{int1}, \eqref{p2} and \eqref{Po1x}   we get
\be\label{a01}
\frac{a_0}2\|u_x(t)\|^2\le \frac k4\|u_{xx}(t)\|^2+\frac{a_0^2}{4k}\|u(t)\|^2\le \left(\frac k4+\frac{a_0h^2}{4k \la_1}\right) \|u_{xx}(t)\|^2 
+\frac{a_0^2h}{4k}\B_N(u(t)).
\ee

Consequently, one has
$$
\E(u(t))\ge \frac m2\|u_t(t)\|^2+\left(\frac k4 -\frac{a_0h^2}{4k\la_1}\right)\|u_{xx}(t)\|^2
+\left(\frac{h\mu}2-\frac{a_0^2h}{4k}\right)\B_N(u(t)).
$$
It is clear that if 
\be\label{muN}
\mu\ge \frac{a^2_0}k \   \mbox{and} \ \ h \  \mbox{is small enough such that} \ \ h^2\le \frac1{2a_0}k^2\la_1,
\ee
then
\be\label{dtE10}
\E(u(t))\ge\frac m2\|u_t(t)\|^2+\frac k8\|u_{xx}(t)\|^2.
\ee
Since the function $a(\cdot)$ satisfies the condition \eqref{ax} then thanks to Poincar\'e inequality
we also have 
\be\label{Ea1}
\E(t)\le\frac m2\|u_t\|^2+\frac12(k+a_1\la_1^{-1})\|u_{xx}\|^2+\frac{\mu h}{2} \B_N(u(t)).
\ee
It follows from \eqref{E1} that
\be\label{215a}
\frac d{dt}\E(u(t))=-b\int_0^L|u_t(x,t)|^{p+2}dx.
\ee
Integrating this equality with respect to $t$ we obtain
$$
\E(u(t))-\E(u(0))+b\int_0^t\int_0^L|u_t(x,s)|^{p+2}dxds =  0.
$$
Hence
\be\label{dtE12}
\int_0^t\int_0^L|u_t(x,s)|^{p+2}dxds\le \frac1b\E(u(0)),
\ee
and 
\be\label{uxxest} 
\frac m2\|u_t(t)\|^2+\frac k8\|u_{xx}(t)\|^2\le \E(u(0)), \ \ \forall t> 0.
\ee
Consider the function 
$$
\E_1(u(t)):=m\|u_t(t)\|^2+\de D_0\|u_{xx}(t)\|^2+\de \int_0^La(x)u^2_x(x,t)dx
+\de h\mu\B_N(u(t)).
$$
Utilising  condition \eqref{ax} and inequality \eqref{p2} we obtain the following inequality
\begin{multline*}
-\de \int_0^La(x)u^2_x(x,t)dx\le \d a_0\|u_x(t)\|^2\\\le \d a_0\|u(t)\| \|u_{xx}(t)\|\le \frac{\d D_0}2 \|u_{xx}(t)\|^2 +\frac{\d a_0^2}{2 D_0} \|u(t)\|^2 \\
\le \frac{\d D_0}2 \|u_{xx}(t)\|^2+\frac{\d a_0^2}{2D_0}\left(h\B_N(u(t)) + \frac{h^2}{\la_1}\|u_{xx}(t)\|^2\right).
\end{multline*}
By using the last inequality we get the estimate
\begin{multline*}
\E_1(u(t))\ge m\|u_t(t)\|^2+\de D_0\|u_{xx}(t)\|^2-\d a_0\|u_x(t)\|^2+\de h\mu\sum_{k=1}^N\overline{u}^2_k(t)\\
\ge m\|u_t(t)\|^2+\left(\frac{\d D_0}2-\frac{\d a_0^2h^2}{2\la_1D_0}\right) \|u_{xx}(t)\|^2+h\left(\de \mu-\frac{\d a_0^2}{2D_0}\right)\B_N(u(t)).
\end{multline*}

If $\mu$  is large enough and $h$ is small enough such that
\be\label{hmL}
h^2\le \frac{\la_1D_0^2}{2a_0^2}, \ \ \mbox{and} \ \ \mu \ge \frac{a_0^2}{D_0},
\ee
then
$$
\E_1(u(t))\ge m\|u_t(t)\|^2+\frac{\de D_0}4\|u_{xx}(t)\|^2+\frac{\mu\de h}{2}\B_N(u(t)).
$$
Furthermore, it follows from the last inequality and  estimate \eqref{Ea1}   that
 \be\label{E1E}  \E_1(u(t))\ge D_1\E(u(t)), \ee
where
$$ 
D_1:=\min\left\{2, \frac{\d D_0}{2(k+a_1\la_1^{-1})}\right\},
$$ 
and according to \eqref{del1} and \eqref{del} 
$\d=\min\left\{\frac{\la_1\sqrt{k}}{4\sqrt{m}}, \ \frac{k\la_1}{\gamma^2}\right\}.$

By using  inequality \eqref{E1E} we obtain from \eqref{dtE9a} that
\begin{multline*}
\frac d{dt}E(u(t))
\le M_0\|u_t(t)\|^2-D_1\E(t)
-b\de\int_0^Lu(x,t)|u_t(x,t)|^pu_t(x,t)dx\\-b\int_0^L|u_t(x,t)|^{p+2}dx,
\end{multline*}
where $M_0:=m(\d+1)+\frac12.$ Integrating last inequality yields
\begin{multline*}
D_1\int_0^t\E(s)ds\le E(u(0))-E(u(t))+M_0\int_0^t\|u_t(s)\|^2ds-b\int_0^t\int_0^L|u(x,s)|^{p+2}dxds\\-b\de\int_0^t\int_0^Lu(x,s)|u_t(x,s)|^pu_t(x,s)dxds.
\end{multline*}
Employing inequality  \eqref{dtE10} and the fact that $\E(t)$ is nonincreasing function we deduce form the last inequality the estimate
\begin{multline}\label{finE}
D_1t\E(u(t))\le E(u(0))+M_0\int_0^t\|u_t(s)\|^2ds-b\int_0^t\int_0^L|u(x,s)|^{p+2}dxds\\-b\de\int_0^t\int_0^Lu(x,s)|u_t(x,s)|^pu_t(x,s)dxds.
\end{multline}
Employing  H\"older's inequality, the Sobolev inequality \eqref{Sob}  and  estimates \eqref{dtE10},  \eqref{dtE12} we get the inequalities
$$
M_0\int_0^t\|u_t(s)\|^2ds\le R_0t^{\frac{2}{p+2}},
$$
where $R_0:=M_0\left[\frac Lb\E(0)\right]^{\frac p{p+2}} ,$
and 
\begin{multline*} 
b\de\left|\int_0^t\int_0^Lu(x,s)|u_t(x,s)|^pu_t(x,s)dxds\right|\\
\le b\d\left(\int_0^t\int_0^L
|u_t(x,s)|^{p+2}dxds\right)^{\frac{p+1}{p+2}}\left(\int_0^t\int_0^L
|u(x,s)|^{p+2}dxds\right)^{\frac1{p+2}}\\
\le b\d\frac{\beta_L}{\la_1}\left(\frac1b\E(0)\right)^{\frac{p+1}{p+2}}\left(\int_0^t\|u_{xx}(s)\|^{p+2}ds\right)^{\frac1{p+2}}\le R_0 t^{\frac1{p+2}},
\end{multline*}
where
$$
R_1:=b\d\frac{\beta_L}{\la_1}\left(\frac1b\E(u(0))\right)^{\frac{p+1}{p+2}}\frac8k\E(u(0)).
$$
Thanks to the last two inequalities we derive from the inequality \eqref{finE}
the desired  decay estimate
$$
D_1E(u(t))\le E(u(0))t^{-1} +R_0 t^{-\frac p{p+2}} +R_1 t^{-\frac{p+1}{p+2}}.
$$
In conclusion we have proved the following theorem:
\begin{theorem} Let $\d=\min\left\{\frac{\la_1\sqrt{k}}{4\sqrt{m}}, \ \frac{k\la_1}{\gamma^2}\right\}$. Suppose  that
 the feedback nudging parameter  $\mu$ is large enough and  $h$ is small enough such that  conditions \eqref{hmu1}, \eqref{muN} and \eqref{hmL} are satsfied, that is, 
$$
h\le \sqrt{\lambda_1}\min\left\{\frac{k}{2a_0},\frac{k}{\sqrt{2a_0}},\frac{D_0}{2a_0}\right\}, \ \  \mu \ge a_0^2\max\left\{\frac1k, \frac1{D_0}\right\},
$$
  then the zero stationary state of  equation \eqref{kg1} is globally asymptotically stable.
Moreover, the following decay estimate is satisfied
$$
\|u_t(t)\|^2+\|u_{xx}(t)\|^2\le C_0 t^{-(p+1)/(p+2)}, \ \ \forall t>0,
$$
where $C_0$ is a positive constant which is independent of  $t.$
\end{theorem}
\section{Stabilization of the  linear equation .}

In this section we are employing  feedback controllers  based on finite-volume elements to stabilize the zero stationary state of the linear wave equation 
$$
mu_{tt}  + ku_{xxxx}  - \left[ {a(x)u_x } \right]_x  + \g u_{tx}
+ bu_t=0,
$$
modeling  dynamics of the  marine riser conveying
fluid (see, e.g., \cite{AdUt} , \cite{LiDeWa} and references therein).
Just as in the previous section, we  assume that the observables are the finite-volume elements $\overline{u}_k(t):=\frac1{|J_k|}\int\limits_{J_k} u(x,t)dx, \ (k=1,\cdots,N)$: 
\be\label{kg1a}
mu_{tt}  + ku_{xxxx}  - \left[ {a(x)u_x } \right]_x  + \g u_{tx}
+ bu_t =
-\mu\sum\limits_{k=1}^N\overline{u}_k(t)\chi_{J_k}(x),    x \in (0,L), \ t>0, 
\ee
\begin{equation}\label{kg2a}
{u}(0,t)={u}(L,t)=u_{x}(0,t)=u_{x}(L,t)=0, \ t>0,
\end{equation}
\begin{equation}\label{kg3a}
u(x,0)=u_0(x), \ \ u_t (x,0)=u_1(x), \ \ x\in (0,L),
\end{equation}
where $k,\g$ and $\mu$ are given positive
parameters.\\
\indent Taking the inner product of 
 the equation  \eqref{kg1} with $u_t+\eb u$ we obtain:

\begin{multline}\label{st3}
\frac d{dt}\Big[\frac m2\|u_t(t)\|^2+\frac k2\|u_{xx}(t)\|^2 +\frac12\int_0^La(x)u_x^2(x,t)dx+\frac12\mu h\B_N(u(t))\\+\eb m(u(t),u_t(t))+\frac{\eb b}2\|u(t)\|^2\Big]
+(b-\eb m)\|u_t\|^2+\eb k\|u_{xx}(t)\|^2\\+\eb\int_0^La(x)u_x^2(x,t)dx+ \mu h \eb\B_N(u(t))=0,
\end{multline} 
where $\B_N(u(t))$ is defined in \eqref{BNu} and  $\eb>0$ is a parameter to be determined below.\\
By using  condition \eqref{ax}, the interpolation inequality \eqref{int1}, the Young's inequality, the inequality \eqref{p2} and the Poincar{\'e} inequality we have
\begin{multline}\label{axin}
-\eb\int_0^La(x)u_x^2(x,t)dx\le \eb a_0\|u_x(t)\|^2\\
 \le  \frac{\eb k}4\|u_{xx}(t)\|^2+\frac{\eb a_0^2}{ k}\|u(t)\|^2\le  \left(\frac{\eb k}4+\frac{\eb a_0^2h^2}{k\la_1}\right)\|u_{xx}(t)\|^2+\frac{\eb a_0^2 h}{ k}\B_N(u(t),
\end{multline}

Employing the  inequality \eqref{axin} we obtain from  \eqref{st3} the inequality

\begin{multline}\label{st4} \frac d{dt} \W(t)+(b-\eb m )\|u_t(t)\|^2
+\left(\frac{3\eb k}4  - \frac{\eb a_0^2h^2}{k\la_1}\right)\|u_{xx}(t)\|^2\\
+\left(\mu h \eb-\frac{\eb a_0^2 h}{ k}\right)\B_N(u(t))\le0,
\end{multline}
where 
\begin{multline}\label{Ygr1}
\W(t):=\frac m2\|u_t(t)\|^2+\frac k2\|u_{xx}(t)\|^2 +\frac12\int_0^La(x)u_x^2(x,t)dx+\frac12\mu h\B_N(u(t))\\+\eb m(u(t),u_t(t))+\frac{\eb b}2\|u(t)\|^2.
\end{multline}

Utilizing the inequalities 
$$
\frac12\int_0^La(x)u_x^2(x,t)dx\ge -\frac{a_0}2\|u_x(t)\|^2
\ge -\frac{a_0}2\|u(t)\|\|u_{xx}(t)\|
\ge -\frac k4\|u_{xx}(t)\|^2-
\frac{a_0^2}{4k}\|u(t)\|^2 ,
$$
and 
$$
\eb m(u(t),u_t(t))\ge -\frac m4 \|u_t(t)\|^2-\frac{\eb^2}m\|u(t)\|^2  ,
$$
 we  we obtain 
$$
\W(t)\ge\frac m4\|u_t(t)\|^2+\frac k4
\|u_{xx}(t)\|^2 +\left(\frac{\eb b}2-\frac {\eb^2}m-\frac{a_0^2}{ 4k}\right)\|u(t)\|^2+\frac12\mu h\B_N(u(t)).
$$
Suppose that
\be\label{emb4}
\eb\le \frac{mb}2.
\ee
Then employing inequality \eqref{p2} yields
$$
\W(t)\ge\frac m4\|u_t(t)\|^2+\left(\frac k4-\frac{a_0^2h^2}{4k\la_1}\right)
\|u_{xx}(t)\|^2 +\left(\frac12\mu h-\frac{a_0^2h}{4k}\right)\B_N(u(t)).
$$

Hence if 
\be\label{muh1}
 h^2\le\frac{k^2\la_1}{2a_0^2} \ \ \mbox{and} \ \ \mu \ge\frac{a_0^2}{k},
\ee
then
\be\label{Fiest2}
\W(t)\ge\frac m4\|u_t(t)\|^2+\frac k8\|u_{xx}(t)\|^2+\frac14\mu h\B_N(u(t)).
\ee
It follows from \eqref{st4} that

\be\label{st4a} \frac d{dt} \W(t)+\frac b2\|u_t(t)\|^2
+\frac{3\eb k}8 \|u_{xx}(t)\|^2
+\frac{a_0^2 h}{\eb k}\B_N(u(t))\le0,
\ee
provided 
$$ \eb \le \frac{b}{2m} , \  h\le \eb k\sqrt{\frac{3\la_1}{8a_0}} \ \mbox{and} \ \ \mu \ge\frac{2a_0^2}{k\eb^2}.
$$
Moreover, it is clear that  there exists a constant $\nu>0$, small enough, such that 
$$
\frac b2\|u_{t}(t)\|^2+\frac{3\eb k}8 \|{u}_{xx}(t)\|^2+
\frac12\mu h\eb \B_N(u(t))\ge \nu \W(t).
$$
Hence one has
$$
\frac d{dt} \W(t)+\nu \W(t)\leq 0.
$$
The above inequality together with  \eqref{Fiest2} imply the exponential
stabilization estimate
\be\label{est2} \|u_t(t)\|^2+\|{u}_{xx}(t)\|^2
\leq D_0 e^{-\nu t}.
\ee
Consequently we have proved the following:
\begin{theorem}\label{T41}
Suppose that 
$$
h\le \min\left\{ \frac{k}{a_0}\sqrt{\frac{\la_1}{2}},  k\sqrt{\frac{3\la_1}{8}}\right\}
\ \ \mbox{and} \ \ 
\mu\ge \frac{a_0^2}{k} \max\left\{\frac12,\frac{2}{\eb^2}\right\},
$$
where $\eb=\min\left\{\frac{mb}{2}, \frac b{2m}.\right\}$
 Then all
solutions of  problem \eqref{kg1a}-\eqref{kg2a} tend to zero with
an exponential rate, as $t \to \infty$.
\end{theorem}

\begin{remark} It is not difficult to see that an analogue of the exponential decay estimate \eqref{est2} is also valid for the marine riser equation with nonliear source term,  i.e., for the equation
$$
mu_{tt}  + ku_{xxxx}  - \left[ {a(x)u_x } \right]_x  + \g u_{tx} 
+ bu_t+f(u) =
-\mu\sum\limits_{k=1}^N\overline{u}_k(t)\chi_{J_k}(x), \  x \in (0,L), \ t>0, 
$$
under the boundary conditions \eqref{kg2a}. Here $f(\cdot)\in C^1(\R)$ is a given nonlinear term that satisfies the conditions 
$$
f(0)=0, \ \ F(s)=\int_0^sf(t)dt\ge0, \ f(s)s-F(s)\ge0, \ \forall s\in \R.
$$

\end{remark}
\begin{remark}  We would like to note that  similarly one can obtain an estimate of the form
$$
\|u_t(t)-v_t(t)\|^2+\|u_{xx}(t)-v_{xx}(t)\|^2\le D_1e^{-\alpha t},
$$
for some $\alpha>0,$ where $u$ is a solution of the following closed-loop feedback control equation

\begin{multline}\label{kg1aa}
mu_{tt}  + ku_{xxxx}  - \left[ {a(x)u_x } \right]_x  + \g u_{tx}
+ bu_t \\=
-\mu\sum_{k=1}^N(\overline{u}_k(t)-\overline{v}_k(t))\chi_{J_k}(x),  x \in (0,L), \ t>0, 
\end{multline}
under the boundary conditions \eqref{kg2a} and initial value \eqref{kg3a}, and $v$ is the solution of the equation
$$
mv_{tt}  + kv_{xxxx}  - \left[ {a(x)v_x } \right]_x  + \g v_{tx}
+ bv_t=0,
$$
under the same boundary conditions \eqref{kg2a}.

\end{remark}


\begin{thebibliography}{9}
\bibitem{Aa}  M. Aasila.  Asymptotic behaviour and stability assessment of marine risers
Math. Methods Appl. Sci. 22 (1999), no. 18, 1585–1598.
\bibitem{AdUt} R. Adiputra and T. Utsunomiya, Stability Based Approach to Design Risers Conveying Fluidfor Ocean Thermal Energy Conversion (OTEC) Application.  Applied Ocean Research,
{\bf 92}, November 2019, 101921
\bibitem{AKN} S. Ahmedov,  V. Kalantarov, A. Namazov. Stabilization of solutions of marine riser equations. Mathematical Methods in the Applied Sciences.  (2024)  https://doi.org/10.1002/mma.10432

\bibitem{AzTi}A.~Azouani, E.S.~Titi, Feedback control of nonlinear dissipative
systems by finite determining parameters - a reaction-diffusion
paradigm,  Evolution Equations and Control Theory, {\bf 3} (2014), 579--594.
\bibitem{AOT} A. Azouani, E. Olson, and E. S. Titi. Continuous data assimilation using general interpolant observables.
J. Nonlinear Sci., {\bf 24(2)} ( 2014) 277–304.
\bibitem{Ba} M.J.~Balas, Feedback control of dissipative hyperbolic distributed parameter systems with finite-dimensional controllers. J. Math. Anal. Appl., {\bf 98(1)} (1984), 1--24.
\bibitem{BaTa} M. Badra and T.  Takahashi, tabilization of parabolic nonlinear systems with finite dimensional feedback or dynamical controllers: Application to the Navier-Stokes system. Journal on Control and Optimization, {\bf 49}(2) (2011) 420–463.
\bibitem{Barbu} V. Barbu,  Boundary stabilization of equilibrium solutions to parabolic equations. IEEE Transactions on Automatic Control, {\bf 58}(9)(2013) 2416–2420.

\bibitem{BaTr} V. Barbu and R. Triggiani,  Internal Stabilization of Navier-Stokes Equations with Finite-Dimensional Controllers, Indiana University
Mathematics Journal,  {\bf  53}(2004) 1443-1494.

\bibitem{Cheb} A. Yu. Chebotarev.  Finite-dimensional controllability of systems of Navier-Stokes type. Differ. Equ. {\bf 46} (2010), no. 10, 1498–1506


\bibitem{CKU1} A.O.Çelebi,  Sh. Gür and V.K. Kalantarov, Structural stability and decay estimate for marine
riser equations. Math. Comput Model. 2011;54: 3182-3
\bibitem{ChLa1} I.D.~Chueshov  and I. Lasiecka, {\it Long-Time Behavior of Second Order Evolution Equations with Nonlinear Damping}, Memories of AMS, no. 912, 2008.
\bibitem{ChLa} I.D.~Chueshov  and I. Lasiecka, {\it Von Karman Evolution Equations: Well-Posedness and Long Time Dynamics},  Springer, New York,2020.

\bibitem{ChKa} I.D.~Chueshov  and V.K.~Kalantarov, Determining functionals for nonlinear damped wave equations. Mat. Fiz. Anal. Geom., {\bf 8(2)} (2001), 215--227.
\bibitem{Coron} J.M. Coron, Control and Nonlinearity, Math.Surveys and Monographs, {\bf 136} AMS, Providence, RI, 2007
\bibitem{CoTr}  J.-M.~ Coron and E.Tr\'{e}lat, Feedback stabilization along a path of steady-states for 1-D semilinear heat and wave equations, Proceedings of the
44th IEEE Conference on Decision and Control, and the European Control Conference 2005
Seville, Spain, December 12-15, 2005.
\bibitem{Far} A. Farhat, E. Lunasin and E.S. Titi, Abridged continuous data assimilation for the 2D Navier–Stokes equations utilizing
measurements of only one component of the velocity field, J. Math. Fluid Mech. 18(1) (2016), 1–23.
\bibitem{FJT} A. Farhat, M.S. Jolly and E.S. Titi, Continuous data assimilation for the 2D Bénard convection through velocity measurements
alone, Phys. D {\bf 303} (2015), 59–66.
\bibitem{Fur} A.V.Fursikov, Optimal Control of Distributed Systems. Theory and Applications, Translations
of Mathematical Monographs, {\bf 187}, Amer. Math. Society, Providence, Rhode Island, 2000
 \bibitem{Gur} S. Gur, Global Asymptotic stability of solutions to nonlinear marine riser equation, Journal of Inequalities and
Applictions,Vol. 2010, Article ID 504670.
).
\bibitem{Kae} S. Kaewunruen
,T. McCarthy,
J. Leklong,
and 
S. Chucheepsakul, Influence of joint stiffness on the free vibrations of a marine riser conveying
fluid. Proceedings of the Eighth (2008) ISOPE Pacific/Asia Offshore Mechanics Symposium
Bangkok, Thailand, November 10-14, 2008 
\bibitem{Ko} M. Köhl, An Extended Liapunov Approach to the
Stabilty Assessment of Marine Risers, Z.Angew. Math. Mech. {\bf 73}
(1993)2,85-92
\bibitem{KaKu} V.K. Kalantarov and A. Kurt,
The Long-Time Behavior of solutions of a Nonlinear Fourth Order Wave
Equation, Describing the Dynamics of Marine Risers, Z. Angew. Math.
Mech. {\bf 77} (1997) 3, 209-215 
\bibitem{KaTi1}  V.~K.~Kalantarov and E.~S.~Titi. Global stabilization of the Navier-Stokes-Voight and the damped nonlinear wave equations by finite number of feedback controllers, Discrete Contin. Dyn. Syst. Ser. B 23 (2018), no. 3, 1325-1345
 \bibitem{KaTi} V.~K.~Kalantarov and E.~S.~Titi. Finite-parameters feedback control for stabilizing damped nonlinear wave equations. \emph{Nonlinear analysis and optimization} 115-133, Contemp. Math., {\bf 659} (2016), Amer. Math. Soc., Providence, RI.
\bibitem{KKV1} H. Kalantarova, V. Kalantarov and O. Vantzos.  Chevron pattern equations: exponential attractor and global stabilization. Vietnam J. Math. {\bf  49} (2021), no. 3, 901–918.
\bibitem{KaOz}  J. Kalantarova, T. Özsarı,. Finite-parameter feedback control for stabilizing the complex Ginzburg-Landau equation.
Systems Control Lett. {\bf 106} (2017), 40–46.
\bibitem{KrSm} M. Krstic and A. Smyshlyaev, Boundary Control of PDEs: A Course on
Backstepping Designs, vol. 16. Philadelphia, PA, USA: SIAM, 2008.

\bibitem{Lar1} A. Larios and Y. Pei. Approximate continuous data assimilation of the 2D Navier–Stokes equations via
the Voigt-regularization with observable data. Evol. Equ. Control Theory, {\bf 9(3)} (2020) 733–751.
\bibitem{LhPr} H. Lhachemi and Chr. Prieur. Finite-dimensional observerbased
boundary stabilization of reaction-diffusion equations with either
a Dirichlet or Neumann boundary measurement. Automatica J. IFAC,
135:Paper No. 109955, 9, 2022.
\bibitem{LiDeWa} Min Li, Di Deng, Decheng Wan, VIV of Flexible Riser Conveying Internal Fluid Subjected to Uniform Current, Proceedings of the Thirtieth (2020) International Ocean and Polar Engineering Conference Shanghai, China, October 11-16, 2020
\bibitem{Lions} J.-L. Lions, {\it Quelques M{\'e}thodes de R{\'e}solution des Probl{\'e}mes aux Limites Non Lin{\'e}aires}.  Dunod, Paris, 1969
\bibitem{LiSt} J.-L. Lions and W. Strauss, Some non-linear evolution equations, Bull. de la Soc. Math. de France, {\bf 93} (1965), 43-96
\bibitem{LiYa} Yun-dong Li,  Yi-ren Yang, Forced vibration of pipe conveying fluid by the Green
function method, Arch Appl Mech (2014) 84:1811–1823
\bibitem{Liu}  W. Liu {\it Elementary Feedback Stabilization of the Linear Reaction-Convection-Diffusion Equation and the Wave Equation}
\bibitem{LuTi} E. Lunasin, E.S. Titi, Finite determining parameters feedback control for distributed nonlinear dissipative systems—a computational study
Evol. Equ. Control Theory 6 (2017), no. 4, 535–557.

\bibitem{Mar}P.A. Markowich, E.S. Titi and S. Trabelsi, Continuous data assimilation for the three-dimensional Brinkman–
Forchheimer-extended Darcy model, Nonlinearity 29(4) (2016), 1292.
\bibitem{Munt} I. Munteanu.{\it Boundary stabilization of parabolic equations}, volume 93
of Progress in Nonlinear Differential Equations and their Applications.
Birkh¨auser/Springer, Cham, 2019. Subseries in Control.
\bibitem{OlTi} E. Olson and E.S. Titi, Determining modes for continuous data assimilation in 2D turbulence, J. Statist. Phys. {\bf 113}(5–
6) (2003), 799–840,
\bibitem{Pai}M.P. Paidosis,  Pipes conveying fluid: A fertile dynamics problem, Journal of Fluids and Structures, {\bf 114} (2022), 2-27. 
\bibitem{Rus}D.L.  Russell, 
Controllability and stabilizability theory for linear partial differential equations: recent progress and open questions
SIAM Review, 20 (4) (1978), pp. 639-739
\bibitem{Sel} A. Selivanov and E. Fridman, Finite-dimensional boundary control of a wave equation with
viscous friction and boundary measurements, IEEE Trans. Automat. Control 69 (2024), no. 5, 3182–
3189.
 \bibitem{Tem} R.~Temam,  \emph{Infinite Dimensional Dynamical Systems in Mechanics and Physics},
New York: Springer, 2nd augmented edition (1997).
\bibitem{Trig}  R. Triggiani, Boundary feedback stabilizability of parabolic equations, Appl. Math. Optim. Applied Mathematics and Optimization, {\bf 6} (1980)201-220
\bibitem{Yan} Y. Yan, D. Coca and V. Barbu, Finite-dimensional controller design for semilinear parabolic systems.
Nonlinear Anal. 70 (2009), no. 12, 4451–4475.
\bibitem{Zua} E. Zuazua {\it Controllability and observability of partial differential equations: some results and open problems} Handbook of differential equations: evolutionary equations, 2011.
 


\end{thebibliography}
\end{document}